\documentclass[12pt]{amsart}
\usepackage{amsmath}
\usepackage{amssymb}

\def\sideremark#1{\ifvmode\leavevmode\fi\vadjust{\vbox to0pt{\vss
  \hbox to 0pt{\hskip\hsize\hskip1em
  \vbox{\hsize2.5cm\tiny\raggedright\pretolerance10000
  \noindent #1\hfill}\hss}\vbox to8pt{\vfil}\vss}}}

\newtheorem{theorem}{Theorem}[section]

\newtheorem{corollary}[theorem]{Corollary}

\theoremstyle{definition}

\newtheorem{question}[theorem]{Question}
\newtheorem{conjecture}[theorem]{Conjecture}
\newtheorem{problem}[theorem]{Problem}

\renewcommand{\H}{\mathbb{H}}

\renewcommand{\phi}{\varphi}
\newcommand{\R}{\mathbb{R}}
\newcommand{\iAH}{{\rm int}(AH(M))}
\newcommand{\rs}{\widehat{\mathbb{C}}}

\begin{document}

\title[Introductory Bumponomics]{Introductory Bumponomics: The topology
of deformation spaces of hyperbolic 3-manifolds}

\author{Richard D. Canary}
\address{University of Michigan \\ Ann Arbor, MI 48109 \\ USA}
\thanks{Research supported in part by grants from the National Science Foundation}

\begin{abstract}
We survey work on the topology of the  space $AH(M)$ of all
(marked) hyperbolic 3-manifolds homotopy equivalent to a fixed
compact 3-manifold $M$ with boundary.
The interior of  $AH(M)$ is quite well-understood, 
but  the topology of the entire space can be quite complicated.
However, the topology is well-behaved at many points in the boundary
of $AH(M)$.
\end{abstract}

\maketitle

\section{Introduction}

In this paper we survey recent work on the topology of the space $AH(M)$
of all hyperbolic 3-manifolds homotopy equivalent to a fixed compact
3-manifold with boundary $M$. The recent resolution of Thurston's Ending Lamination
Conjecture (in Minsky \cite{ELC1} and Brock-Canary-Minsky \cite{ELC2,ELC3})  in combination
with the resolution of Marden's Tameness Conjecture (in Agol \cite{agol}
and Calegari-Gabai \cite{calegari-gabai}),
gives a complete classification of the manifolds in $AH(M)$.
However, the invariants in this classification vary discontinuously
over the space, and we are very far from having a parameterization of $AH(M)$.

The interior $\iAH$ has been well-understood since the 1970's,
due to work of
Ahlfors, Bers, Kra, Marden, Maskit, Sullivan, Thurston and others. The components of
$\iAH$ are enumerated by topological data, and each component is a manifold parameterized
by natural analytic data. The recent resolution of the Bers-Sullivan-Thurston
Density Conjecture assures us that $AH(M)$ is the closure of $\iAH$.

Since the mid-1990's, there has been a string of results and examples
demonstrating that the topology of $AH(M)$ is less well-behaved than originally
expected. Anderson and Canary \cite{ACpages} first showed that components
of $\iAH$ can bump, i.e. have intersecting closure,
while Anderson, Canary and McCullough \cite{ACM}
characterized exactly which components of $\iAH$ can bump when $M$ has 
incompressible boundary. McMullen \cite{mcmullen} showed that
if $M=S\times I$ (where $S$ is a closed surface), then the only
component of $\iAH$ self-bumps, i.e. there is a point in the boundary
such that any sufficiently small neighborhood of the point disconnects the 
interior. Bromberg and Holt \cite{bromberg-holt} showed that self-bumping
occurs whenever $M$ contains a primitive essential annulus.
Most recently, Bromberg \cite{bromberg-PT} showed that the space of
punctured torus groups is not even locally connected.
We will survey these results and describe the construction which
has been responsible for all the pathological behavior discovered so far.

On the other hand, the topology of the deformation space appears to be well-behaved
at most points on its boundary and we will describe some recent results establishing
this in a variety of settings.

{\bf Acknowledgements:} I would like to thank the organizers of the workshop, Ravi Kulkarni
and Sudeb Mitra, for
inviting me to visit the Harish-Chandra Research Institute. I enjoyed both the 
mathematical interaction and the opportunity to visit India for the first time. I would especially like
to thank the student organizers, Vikram Aithal,
Krishnendu Gongopadhyay, and 
Siddhartha Sarkar, for their help during my stay.

I would also like to thank the referee for their very careful reading of the original
version of this manuscript and their many useful suggestions.

\section{Definitions}

In this section, we will set up some of the notation and introduce some of
the definitions used throughout the remainder of the paper.

Let $M$ be a compact 3-manifold whose interior
admits a complete hyperbolic structure.  We will assume throughout this paper
that all surfaces and 3-manifolds are oriented and have non-abelian fundamental
group and that all manifolds are allowed to have boundary. Then
$AH(M)$ is the space of (marked) hyperbolic 3-manifolds homotopy equivalent
to $M$. More formally, we define 
$$AH(M) = \{\rho:\pi_1(M)\to {\rm PSL}_2({\bf C})|\ \  \rho\ 
{\rm discrete}\ {\rm and\ faithful }\}/{\rm PSL}_2({\bf C}).$$
We topologize $AH(M)$ as a subset of the the character variety
 $$X(M)=Hom_T(\pi_1(M),{\rm PSL}_2({\bf C}))//{\rm PSL}_2({\bf C})$$
 where $Hom_T(\pi_1(M),{\rm PSL}_2({\bf C}))$ denotes the space of
representations  $\rho:\pi_1(M)\to{\rm PSL}_2({\bf C})$ with
 the property that if $g$ is an element of a rank two abelian subgroup of
 $\pi_1(M)$, then $\rho(g)$ is either parabolic or the identity.
(Here, we are taking the Mumford quotient to guarantee that the quotient
has the structure of an algebraic variety, see Kapovich \cite{kapovich-book} for details.)

If $\rho\in AH(M)$, then $N_\rho=\H^3/\rho(\pi_1(M))$ is a hyperbolic 3-manifold
homotopy equivalent to $M$. (The elements of $AH(M)$ are really equivalence
classes of representations, but as the manifolds associated to conjugate
representations are isometric we will consistently blur this distinction.)
There is also a homotopy equivalence $h_\rho:M\to N_\rho$ such that
$$ (h_\rho)_*:\pi_1(M)\to \pi_1(N_\rho)=\rho(\pi_1(M))$$
agrees with $\rho$. The homotopy equivalence $h_\rho$ is the ``marking'' of $N_\rho$.

Alternatively, we could have defined $AH(M)$ as the space of
pairs $(N,h)$ where $N$ is an oriented  hyperbolic 3-manifold and $h:M\to N$ is
a homotopy equivalence, where we consider two pairs $(N_1,h_1)$ and
$(N_2,h_2)$ to be equivalent if there is an orientation-preserving isometry
$j:N_1\to N_2$ such that $j\circ h_1$ is homotopic to $h_2$. This alternate definition
is reminiscent of the classical definition of Teichm\"uller space as a space of marked
Riemann surfaces of a fixed genus.

If $\rho \in AH(M)$, then the {\em domain of discontinuity} $\Omega(\rho)$ is the  largest
open subset of $\hat {\bf C}$ on which
$\rho(\pi_1(M))$ acts properly discontinuously. 
The {\em limit set} $\Lambda(\rho) =\hat{\bf C} -\Omega(\rho)$ is its complement.
The {\em conformal boundary} is defined to be
$\partial_cN_\rho=\Omega(\rho)/\rho(\pi_1(M))$ and we let
$$\hat N_\rho=N_\rho\cup \partial_cN_\rho=(\H^3\cup \Omega(\rho))/\rho(\pi_1(M)).$$

We say that $N_\rho$ is {\em convex cocompact}, if $\hat N_\rho$ is compact. We say
that $N_\rho$ is {\em geometrically finite} if $\hat N_\rho$ is homeomorphic to
$M'-P'$ where $M'$ is a compact 3-manifold and $P'$ is a collection of annuli and
tori in $\partial M'$. (These definitions are equivalent to more classical definitions,
see Marden \cite{marden} and Bowditch \cite{bowditch}.)

\section{The interior of $AH(M)$}
\label{interior}

In this section, we survey the classical deformation theory of hyperbolic
3-manifolds which gives a beautiful description of the interior $\iAH$ of $AH(M)$. If the boundary of $M$ consists
entirely of tori, then Mostow-Prasad Rigidity \cite{mostow,prasad} implies that any homotopy
equivalence between hyperbolic 3-manifolds homotopy equivalent to $M$ is homotopic
to an isometry. Therefore, $AH(M)$ has either
0 or 2 points (one gets one point for each orientation of $M$.) So, we will always assume that the boundary of $M$ has a non-toroidal component.

Marden \cite{marden} and Sullivan \cite{sullivan2} proved that $\iAH$ consists exactly of the
geometrically finite hyperbolic manifolds $\rho\in AH(M)$ such that every parabolic
element of $\rho(\pi_1(M))$ is contained in a free abelian subgroup of rank two.
(This is equivalent to requiring that $\hat N_\rho$ be homeomorphic to a compact
3-manifold with its toroidal boundary components removed.) 
The now classical quasiconformal deformation theory of Kleinian groups shows
that geometrically finite hyperbolic 3-manifolds are
determined by their (marked) homeomorphism type and  the conformal
structure on the conformal boundary and that every possible conformal
structure arises. For an analytically-oriented discussion of quasiconformal
deformation theory, see Bers \cite{BersSurvey}. For a more topological
viewpoint on this material, see chapter 7 of Canary-McCullough \cite{canary-mccullough}.

In order to formally state the parameterization theorem for $\iAH$ we need
to introduce some more notation.

We define ${\mathcal A}(M)$ to be the set of oriented, compact, irreducible,
atoroidal (marked) 3-manifolds homotopy equivalent to $M$. More formally,
$\mathcal{A}(M)$ is the set of pairs $(M',h')$ where $M'$ is an oriented,
compact, irreducible, atoroidal 3-manifold and $h':M\to M'$ is a homotopy
equivalence where two pairs $(M_1,h_1)$ and $(M_2,h_2)$ are considered
equivalent if there exists an orientation-preserving homeomorphism $j:M_1\to M_2$
such that $j\circ h_1$ is homotopic to $h_2$. We recall that $M'$ is said to be
{\em irreducible} if every embedded 2-sphere in $M'$ bounds a ball and is 
said to  be {\em atoroidal} if every rank two free abelian subgroup of $\pi_1(M')$ is
conjugate to  a subgroup of $\pi_1(T')$ for some toroidal boundary component
$T'$ of $M'$.

If $[(M',h')]\in \mathcal{A}(M)$, we define $Mod_0(M')$ to be the group of isotopy
classes of homeomorphisms of $M'$ which are homotopic to the identity. We define
$\partial_{NT}M'$ to be the non-toroidal components of $\partial M$ and we let
${\mathcal T}(\partial_{NT}M')$ denote the Teichm\"uller space of all (marked) conformal
structures on $\partial_{NT}M'$. Recall that the Teichm\"uller space of a disconnected
surface is simply the product of the Teichm\"uller spaces of its components, so
${\mathcal T}(\partial_{NT}M')$ is always topologically a cell.

\begin{theorem}{}{(Ahlfors, Bers, Kra, Marden, Maskit, Sullivan,Thurston)}{}
$$\iAH \cong \bigcup_{[(M',h')]\in {\mathcal A}(M)} {\mathcal T}(\partial_{NT}M')/Mod_0(M')$$
\end{theorem}

This identification is quite natural. By the previously mentioned results of 
Marden and Sullivan, if $\rho \in \iAH$, there exists a compact, atoroidal, irreducible
3-manifold $M_\rho$ and an orientation-preserving homeomorphism 
$j_\rho:\hat N_\rho\to {\rm int}(M_\rho)\cup \partial_{NT}M_\rho$, so we obtain
a well-defined marked homeomorphism type $[(M_\rho,j_\rho\circ h_\rho)]\in {\mathcal A}(M)$.
If \hbox{$[(M_\rho,j_\rho\circ h_\rho)]=[(M',h')]$,} then we may assume that $M_\rho=M'$ and choose $j_\rho$ so that
$j_\rho\circ h_\rho$ is homotopic to $h'$. With this constraint, $j_\rho$ is well-defined up to post-composition by elements
of $Mod_0(M')$ and
$\partial_cN_\rho$ is a Riemann surface,
so we get  a well-defined element of $ {\mathcal T}(\partial_{NT}M')/Mod_0(M')$.

Let $\Phi: \iAH\to \bigcup_{[(M',h')]\in {\mathcal A}(M)} {\mathcal T}(\partial_{NT}M')/Mod_0(M')$ be the map
defined in the previous paragraph. We will discuss the ingredients of the proof that $\Phi$ is a homeomorphism.
Maskit's Extension Theorem \cite{maskit} implies that if
$\Theta(\rho_1)=\Theta(\rho_2)$, then there exists a quasiconformal map $\phi:\rs\to\rs$ which is
conformal on $\Omega(\rho_1)$ and conjugates the action of $\rho_1(\pi_1(M))$ to the action of
$\rho_2(\pi_1(M))$. Ahlfors \cite{AhlforsMZ} proved that if $\rho_1$ is geometrically finite, then
$\Lambda(\rho_1)$ has measure zero, so we may conclude  that  $\phi$ is conformal. Therefore,
$\Phi$ is injective. 
The surjectivity of $\Phi$ follows from the Measurable
Riemann Mapping Theorem  \cite{ABMRT} and Thurston's Geometrization Theorem
(see Morgan \cite{Morgan}). 

If $M$ has incompressible boundary (equivalently if
$\pi_1(M)$ is freely indecomposable), then $Mod_0(M')$ is trivial, so each component of
${\rm int}(AH(M))$ is topologically a ball.
Moreover,  Canary and McCullough \cite{canary-mccullough}, showed that  if
$M$ has incompressible boundary, then $\mathcal{A}(M)$ is infinite if and only if
$M$ has {\em double trouble}, i.e. there exist  simple closed curves $\alpha$ and $\beta$ in
$ \partial_{NT}M$ which are both homotopic to  a  curve $\gamma$ in a toroidal boundary component
of $M$ but are not homotopic in $\partial M$.  (Alternatively, $M$ has double trouble if and only if there
is a  thickened torus component of its characteristic submanifold, which intersects
the boundary of $M$ in at least two annuli.)
We summarize this below.

\begin{corollary}{}{} If $M$ has incompressible boundary, then $\iAH$ is
homeomorphic to a collection of disjoint balls. This collection is infinite if and only
if $M$ has double trouble.
\end{corollary}

If $M$ has compressible boundary, then typically
$\mathcal{A}(M)$ is infinite (see Canary-McCullough \cite{canary-mccullough} for a detailed
analysis) and $Mod_0(M)$ is infinitely generated (see McCullough \cite{mccullough-twist}.)
Maskit \cite{maskit} showed that $Mod_0(M')$ always acts freely on $\mathcal{T}(\partial_{NT}M')$, so each
component of $\iAH$ is still a manifold.

\medskip\noindent
{\bf Examples: 1)} If $M$ is homeomorphic to a trivial $I$-bundle $S\times I$ over a closed
orientable surface $S$ of genus $g\ge 2$, then ${\mathcal A}(M)$ has only one element.
(Notice that $S\times I$ has the quite rare property that it admits an orientation-reversing
self-homeomorphism homotopic to the identity.) The elements of $\iAH$ are called
{\em quasifuchisan} as they are quasiconformally conjugate to Fuchsian groups.
In this case, $\iAH$ is often denoted $QF(S)$ and $QF(S)\cong {\mathcal T}(S)\times{\mathcal T}(S)$. Concretely, if
$\rho\in\iAH$, then $N_\rho$ is determined by the conformal structure on
$\partial_cN_\rho$ and given any conformal structure on $S\times \{ 0,1\}$ one
can construct  a complete hyperbolic structure on $S\times (0,1)$ with appropriate
conformal boundary.

\medskip
\noindent
{\bf 2)} Let $S_i$ denote a surface of genus $i$ with one (open) disk removed.
We consider the $I$-bundle $J_i=S_i\times [0,1]$ and let $\partial_r J_i=\partial S_i\times [0,1]$.
Let $\{ A_i\}_{i=1}^n$ be a collection of disjoint, parallel, consecutively ordered
longitudinal annuli in the
boundary of a solid torus $V$. Then $M_n$ is formed from $V$ and $\{ J_1,\ldots, J_n\}$
by attaching $\partial_rJ_i$ to $A_i$. The manifolds $M_n$ are examples of
books of $I$-bundles, see Culler-Shalen \cite{culler-shalen}.

Any irreducible manifold homotopy equivalent
to $M_n$ is formed by attaching the $\{ J_i\}$ to the $\{ A_i\}$ in a different order.
To be more precise, if $\tau$ lies in the permutation group $\Sigma_n$, then one
can form $M^\tau_n$ from $V$ and $\{ J_1,\ldots, J_n\}$ by attaching
$\partial_rJ_i$ to $A_{\tau(i)}$. One may extend the identity map on $\{ J_1,\ldots, J_n\}$
to a homotopy equivalence $h_\tau:M_n\to M^\tau_n$. It is a consequence, see \cite{ACpages}, of
Johannson's  Deformation Theorem \cite{johannson},  that every
$(M',h')\in \mathcal{A}(M_n)$ is equivalent to $(M^\tau_n,h_\tau)$ for some $\tau\in\Sigma_n$.

It is easy to check that if $n\ge 3$, $(M^\tau_n,h_\tau)$ need not
be equivalent to $(M_n,id)$. For example, if $n=4$ and $\tau=(2\ 3)$ one may check that
the boundary components of $M_4$ have genera 3, 5, 5, and 7, while the boundary components
of $M^\tau_4$ have genera 4, 5, 5, and 6, so $M_4$ and $M^\tau_4$ are clearly not
homeomorphic. On the other hand, if $\tau$ is any multiple of the rotation $(1\ 2\cdots n)$,
$(M^\tau_n,h_\tau)$ is clearly equivalent to $(M_n,id)$. Analyzing the situation more 
carefully, again using Johannson's Deformation Theorem, one shows that
the space $\mathcal{A}(M_n)$  of
marked homeomorphism types of manifolds homotopy equivalent to $M_n$ may be identified
with cosets of  the subgroup generated by $(1\ 2\cdots n)$ in the permutation group $\Sigma_n$
(see \cite{ACpages}). In particular, $\mathcal{A}(M_n)$
has $(n-1)!$ elements, so ${\rm int}(AH(M_n))$ is homeomorphic to a disjoint union of
$(n-1)!$ balls. Notice that $M_2=S_3\times I$ where $S_3$ is a closed orientable surface
of genus 3.

\section{The Density Theorem}

Much of the recent work on infinite volume hyperbolic 3-manifolds has
culminated in the proof of the Bers-Sullivan-Thurston Density
Conjecture which asserts that $AH(M)$ is the closure of its
interior. More concretely, the Density Conjecture predicts that
every hyperbolic 3-manifold with finitely generated fundamental group
is an (algebraic) limit of geometrically finite hyperbolic 3-manifolds.

\medskip\noindent
{\bf Density Theorem:} {\em  If $M$ is a compact hyperbolizable
3-manifold, then $AH(M)$ is the closure of its interior $\iAH$.}

\medskip

There are two approaches to the proof of the Density Theorem.
Both approaches make use of the proof of Marden's Tameness
Conjecture by Agol \cite{agol} and Calegari-Gabai \cite{calegari-gabai}
which asserts that every hyperbolic 3-manifold with finitely generated
fundamental group is homeomorphic to the interior of a compact 3-manifold.

In the classical approach one  starts with a hyperbolic manifold $N_\rho\in AH(M)$
and constructs a sequence $\{ N_{\rho_i}\}$ 
in $\iAH$ such that the end invariants of $\{N_{\rho_i}\}$ converge to
those of $N_\rho$. One then
uses convergence results of Thurston \cite{thurston-notes},
Kleineidam-Souto\cite{kleineidam-souto}, Lecuire \cite{lecuire} and
Kim-Lecuire-Ohshika \cite{KLO}  to  show that this sequence converges to
a limit $N_{\rho'}$. Arguments of Namazi-Souto \cite{namazi-souto} or
Ohshika \cite{ohshika-density} can then be used to check that $N_{\rho'}$ has the
same end invariants as $N_\rho$. One then invokes the solution of
Thurston's Ending Lamination Conjecture \cite{ELC1,ELC2,ELC3} to show that $N_\rho=N_{\rho'}$
and thus complete the proof that $N_\rho\in\overline{\iAH}$.

The other approach makes use of the deformation theory of cone-manifolds
developed by Hodgson-Kerckhoff \cite{HK1,hodgson-kerckhoff} 
and Bromberg \cite{bromberg-deform}. These ideas were first used by Bromberg
\cite{bromberg-density} to prove Bers' original Density conjecture for
hyperbolic manifolds without cusps. This conjecture asserted that any Kleinian
surface group with exactly one invariant component of its domain of discontinuity
was in the boundary of the ``appropriate'' Bers slice (see section \ref{relative} for the
definition of a Bers slice). Brock and Bromberg 
\cite{brock-bromberg} showed that if $M$ has incompressible boundary and
$N_\rho\in AH(M)$ has no cusps, then $N_\rho\in \overline{\iAH}$. A complete
proof of the Density Theorem using these methods was given by Bromberg
and Souto \cite{bromberg-souto}.

\section{Bumping}

It would be reasonable to expect that no two components of ${\rm int}(AH(M))$
have intersecting closures, since one might expect that any hyperbolic manifold
in the closure of a component $B$ would be homeomorphic to a hyperbolic manifold
in $B$. In fact, Jim Anderson and I spent two years attempting to prove this
and eventually came up with examples \cite{ACpages} which illustrated that this ``bumping''
of components can occur. Later, Anderson, Canary and McCullough \cite{ACM}
characterized exactly when two components of ${\rm int}(AH(M))$ can bump in the
case that $M$ has incompressible boundary.

Formally, 
we say that two components $B$ and $C$ of $\iAH$ {\em bump at $\rho$} if
$\rho\in \overline{B}\cap \overline{C}$. Notice that whenever two components of $\iAH$ bump,
$AH(M)$ is not a manifold.

The first setting in which the phenomenon of bumping was observed, was the books $\{ M_n\}$ of $I$-bundles
discussed in section \ref{interior}.

\begin{theorem}{\rm (Anderson-Canary \cite{ACpages})} 
\label{pages} If $n\ge 3$, then any two components of 
${\rm int}(AH(M_n))$ bump, where $\{ M_n\}$ are the books of I-bundles constructed in 
example 2 in section \ref{interior}. In particular, $AH(M_n)$ is connected, while ${\rm int}(AH(M_n))$
has $(n-1)!$ components.
\end{theorem}

We will outline the construction in the proof of Theorem \ref{pages} in the section \ref{key}. 
In a finer analysis,  Holt \cite{holtCAG} showed that there is a single point at which any two components  of  ${\rm int}(AH(M_n))$ bump.

\begin{theorem}{\rm (Holt \cite{holtCAG})} If $n\ge 3$, then there exists $\rho\in AH(M_n)$
which  lies in the boundary of every component of ${\rm int}(AH(M_n))$.
\end{theorem}

{\bf Remark:} Holt \cite{holtCAG} further observes, that the set of points lying in the closure
of every component  of ${\rm int}(AH(M_n))$ contains an open subset of a (complex) codimension 1 subvariety of $X(M_n)$. It is typical that the ``bumping locus'' is relatively large.

\medskip

Anderson, Canary and McCullough \cite{ACM} later gave a complete characterization of
which components of $\iAH$ can bump in the case when $M$ has incompressible
boundary. Roughly, the two components can bump if the two (marked) homeomorphism
types differ by cutting along a collection of (primitive) solid tori and rearranging
the order in which the complementary pieces are attached. One sees that for $M_n$
any two homeomorphism types differ in this specific way.
A solid torus $V$ in $M$ is said to be {\em primitive} if $V\cap \partial M$ is a  non-empty
collection of annuli, the inclusion of each annulus of $V\cap \partial M$ 
into $V$ is a homotopy equivalence and the image of $\pi_1(V)$ in  $\pi_1(M)$
is a maximal abelian subgroup. 

To illustrate the role of primitivity in Anderson, Canary and McCullough's 
result, we consider a sequence of manifolds $\{ M_n'\}_{n=3}^\infty$,
again obtained from a solid torus $V$ and
$I$-bundles $\{ J_1,\ldots,J_n\}$. This time we let $\{ A_1',\ldots,A_n'\}$ be a collection of
disjoint, parallel, consecutively ordered annuli in the boundary of a solid torus $V$ such
that the inclusion of $\pi_1(A_i')$ into $\pi_1(V)$ is a subgroup of index 3 (i.e. each annulus
wrap 3 times around the longitude of $V$.) We form $M_n'$ by attaching $\partial_rJ_i$
to $A_i'$. It is again the case that any manifold homotopy equivalent to $M_n'$ is obtained
by attaching the $J_i$ in a different order and that $\mathcal{A}(M_n')$ has
$(n-1)!$ elements. However, the results of \cite{ACM} imply that no two components
of ${\rm int}(AH(M_n'))$ bump for any $n$.

In order to give a precise statement of the results of \cite{ACM} we must
introduce the notion of a  primitive shuffle equivalence. In what follows, if
$M$ is a compact, irreducible 3-manifold with incompressible boundary, then
$\Sigma(M)$ will denote its characteristic submanifold. For complete discussions
of the theory of characteristic submanifolds see Jaco-Shalen \cite{jaco-shalen} or
Johannson \cite{johannson}. For a discussion in the setting of hyperbolizable 3-manifolds,
see Canary-McCullough \cite{canary-mccullough} or Morgan \cite{Morgan}.

Given two compact irreducible $3$-manifolds $M_1$ and $M_2$
with nonempty incompressible boundary, a homotopy equivalence
$h\colon\, M_1\to M_2$ is a {\em
primitive shuffle equivalence} if there exists a finite collection $V_1$ of
primitive solid torus components of $\Sigma(M_1)$ and a finite
collection $V_2$ of solid torus components of $\Sigma(M_2)$, so
that $h^{-1}(V_2)=V_1$ and so that $h$ restricts to an
orientation-preserving homeomorphism from the closure of
$M_1-V_1$ to the closure of $M_2-V_2$.

If $M$ is a compact, hyperbolizable 3-manifold with nonempty
incompressible boundary, we say that two elements $[(M_1,h_1)]$ and
$[(M_2,h_2)]$ of $\mathcal{A}(M)$ are {\em primitive shuffle equivalent}
if there exists a primitive shuffle equivalence $s\colon\, M_1\to M_2$ such that
$[(M_2,h_2)] =[(M_2, s\circ h_1)]$. 

\begin{theorem}{}{\rm (Anderson-Canary-McCullough \cite{ACM})}{}
\label{PSE}
Let $M$ be a compact, hyperbolizable
$3$-manifold with nonempty incompressible boundary, and let
$[(M_1,h_1)]$ and $[(M_2,h_2)]$ be two elements of $\mathcal{ A}(M)$. The
associated components of $\iAH$ have intersecting closures if
and only if $[(M_2,h_2)]$ is primitive shuffle equivalent to
$[(M_1,h_1)]$.
\end{theorem}

Combining the work of Anderson, Canary and McCullough \cite{ACM} with the
resolution of the Density Conjecture, one obtains
a complete enumeration of the components of $AH(M)$ when $M$ has
incompressible boundary. In particular, one completely determines
exactly when $AH(M)$ has infinitely many components.  
Primitive shuffle equivalence gives a finite-to-one
equivalence relation on $\mathcal{A}(M)$ and we let $\widehat{\mathcal{A}}(M)$
be the quotient of $\mathcal{A}(M)$ by this equivalence relation.

\begin{corollary}{}{}
\label{enum}
If $M$ has incompressible boundary, then the components of $AH(M)$ are in one-to-one
correspondence with $\widehat{\mathcal{A}}(M)$. In particular, $AH(M)$ has infinitely many
components if and only if $M$ has double trouble.
\end{corollary}

Holt \cite{holtAJM} refined the analysis of \cite{ACM} to show that if $C_i$ is a collection of components of $\iAH$ such that any two components in the collection bump, then
then they all bump at a single point. 

\begin{theorem}{\rm (Holt \cite{holtAJM})}
Let $M$ be a compact, hyperbolizable
$3$-manifold with nonempty incompressible boundary, and let
$\{[(M_i,h_i)]\}_{i=1}^m$ be a collection of elements of $\mathcal{ A}(M)$ such that
$[(M_i,h_i)]$ is primitive shuffle equivalent to $[(M_1,h_1)]$ for all $i$. If $C_i$ is
the component of $\iAH$ associated to $[(M_i,h_i)]$,  then
$$\bigcap_{i=1}^m \overline{C_i}\ne\emptyset.$$
\end{theorem}

\section{Self bumping}

McMullen \cite{mcmullen} was the first to observe that individual components of
$\iAH$ can self-bump.  
A component $B$ of ${\rm int}(AH(M))$ {\em self-bumps} at $\rho\in \overline{B}$
if there is a neighborhood $V$ of $\rho$  such that if $\rho\in W\subset V$ is any
sub-neighborhood, then $W\cap B$ is disconnected.

\begin{theorem}{\rm (McMullen \cite{mcmullen})} If $S$ is a closed surface,
then $QF(S)={\rm int}(AH(S\times I))$ self-bumps.
\end{theorem}

Notice that this implies, in particular, that $AH(S\times I)$ is not
a manifold. The self-bumping points in \cite{mcmullen} are again obtained using the construction
from \cite{ACpages} described in section \ref{key}. 

McMullen's proof made use, in a crucial manner, of the theory of projective structures on
surfaces, so did not generalize immediately to manifolds which are not I-bundles. Bromberg
and Holt \cite{bromberg-holt} were able to generalize McMullen's result to all manifolds
which contain primitive essential annuli. An embedded annulus $A$ in $M$ is said
to be {\em essential} if  it is {\em incompressible}, i.e. $\pi_1(A)$ injects into $\pi_1(M)$,
and is not properly homotopic into the boundary $\partial M$. It is said to be {\em primitive}
if the image of $\pi_1(A)$ in $\pi_1(M)$ is a maximal abelian subgroup of $\pi_1(M)$.

\begin{theorem}{\rm (Bromberg-Holt \cite{bromberg-holt})}
\label{BH}
If $M$ contains a primitive essential
annulus, then every
component of $\iAH$ self-bumps.
\end{theorem}

Notice that there is no assumption that $M$ has incompressible boundary in this
theorem. It implies, in particular, that $AH(M)$ is not a manifold if $M$ contains a primitive
essential annulus. 

\medskip

Ito has completed an extensive analysis of related phenomena in the space $P(S)$ of
complex projective structures on a surface $S$. There is a natural map $hol:P(S)\to X(S\times I)$
which takes a complex projective structure to its associated holonomy map. Hejhal \cite{hejhal}
showed that the map $hol$ is a local homeomorphism and 
Goldman \cite{goldman} showed that the components of $Q(S)=hol^{-1}(QF(S))$ are enumerated
by the set of weighted multicurves on $S$. McMullen \cite{mcmullen} showed that
$QF(S)$ self-bumps by showing that two components of $Q(S)$ can bump. Ito \cite{ito-proj1,ito-proj2}
shows that any two components of $Q(S)$ bump, that any component of $Q(S)$ other than
the base component self-bumps, and that arbitrarily many components of $Q(S)$ can bump
at a single point. (Bromberg and Holt have obtained related results.)

\section{The key construction}
\label{key}

In this section, we will describe the ``wrapping'' construction of
\cite{ACpages} in
the special case of $AH(M_4)$ and show that components of ${\rm int}(AH(M_4))$ can bump.
 
Let $\tau\in\Sigma_4$ be the permutation $(2\ 3)$ and
let $\hat M^\tau_4$ be the manifold obtained from $M^\tau_4$ by removing an open
neighborhood of the core curve of $V$. One may construct an infinite cyclic cover 
$(\hat M^\tau_4)'$ of $\hat M^\tau_4$
from an infinite thickened annulus $S^1\times I\times \R$ by attaching
infinitely many copies of each $J_i$ to the outer boundary
$S^1\times \{0\}\times \R$ so that they occur repeatedly in the cyclic order
$\ldots J_1,J_3,J_2,J_4,J_1,J_3,\ldots...$. (More concretely, for all
$n\in {\bf Z}$ and $i=1,2,3,4$ one attaches
a copy of $J_i$ to the thickened annulus by identifying $\partial_rJ_i$ with
$S^1 \times \{ 0\} \times [12n +3\tau(i)-1,12n+3\tau(i)+1]$ by an
orientation-reversing homeomorphism. 
Vertical translation by 12 on the infinite thickened annulus extends to a
homeomorphism of 
$(\hat M^\tau_4)'$ which generates the full group of covering transformations of $(\hat M^\tau_4)'$ over
$M^\tau_4$.) Let $\pi:(\hat M^\tau_4)'\to \hat M^\tau_4$ be the covering map.

One then constructs an orientation-preserving embedding $\hat f_\tau: M_4\to(\hat M^\tau_4)'$  which takes each
$J_i$ homeomorphically to a copy of $J_i$. (More concretely, one may take $J_i$ to the copy
of $J_i$ attached to $S^1 \times I \times [12i +3\tau(i)-1,12i+3\tau(i)+1]$.) Let $f_\tau=\pi\circ\hat f_\tau$.

Let $\tilde M^\tau_4$ be the cover of $\hat M^\tau_4$ associated to $(f_\tau)_*(\pi_1(M_4))$.
One easily checks that $f_\tau$ lifts to an embedding $\tilde f_\tau:M_4\to\tilde M^\tau_4$
(since it lifts to an embedding in the intermediate cover $(\hat M^\tau_4)'$.) Also, notice
that if $i_0:\hat M^\tau_4\to M^\tau_4$ denotes the inclusion map, then
$i_0\circ f_\tau$ is a homotopy equivalence and is homotopic to $h_\tau$
 (where $h_\tau$ is the homotopy equivalence
defined in example 2 in section \ref{interior}).

The central tool in the construction is the generalization of Thurston's
Hyperbolic  Dehn Filling
Theorem \cite{thurston-notes} to the setting of geometrically finite hyperbolic 3-manifolds.

Let $T$ be a toroidal boundary component of  compact 3-manifold $M$
and let $(m,l)$ be  a choice of  meridian and longitude for $T$.
Given a pair $(p,q)$ of relatively
prime integers, we may form a new manifold $M(p,q)$ by attaching a solid torus $V$ 
to $M$ by an orientation-reversing homeomorphism $g\colon\, \partial V \to T$
so that, if $c$ is the meridian of $V$, then $g(c)$ is a $(p,q)$
curve on $T$ with respect to the chosen meridian-longitude system.
We say that $M(p,q)$ is obtained from $M$ by {\em
$(p,q)$-Dehn filling along $T$}.

If $T=\partial M$ and $M$ is hyperbolizable (i.e. its interior admits a complete
hyperbolic structure), then Thurston \cite{thurston-notes} proved that 
$M(p,q)$ is hyperbolizable for all but finitely many choices of $(p,q)$.
Bonahon and Otal \cite{bonahon-otal} were the first to observe that you
could generalize this result to the setting of geometrically finite hyperbolic
3-manifolds (see also Comar \cite{comar}, Hodgson-Kerckhoff \cite{HK1} and
Bromberg \cite{bromberg-ds}.) Our statement of the resulting Hyperbolic Dehn
Filling Theorem is rather convoluted but it essentially says that given
a geometrically finite hyperbolic 3-manifold $N$ homeomorphic to ${\rm int}(M)$ and
$(p,q)$ near to $\infty$, then there exists a geometrically finite hyperbolic
3-manifold $N(p,q)$ which is ``near'' to $N$.
The complications in the statement are largely the result of the need
to keep careful track of the marking.

\medskip\noindent
{\bf Hyperbolic Dehn Filling Theorem:} {\em
Let $M$ be a compact, hyperbolizable $3$-manifold and let $T$ be toroidal
boundary component of $M$.
Let $N=\H^3/\Gamma$ be a geometrically finite
hyperbolic $3$-manifold admitting an orientation-preserving homeomorphism
 $\psi:{\rm int}(M) \to  N$. Further assume that every
parabolic element of $\Gamma$ lies in a rank-two parabolic
subgroup. Let $\{(p_n,q_n)\}$ be an infinite sequence of distinct relatively
prime pairs of integers.

Then, for all sufficiently large $n$, there exists a (non-faithful)  representation
$\beta_n:\Gamma\to {\rm PSL}_ 2({\bf C})$ with discrete image such that
\begin{enumerate}
\item $\beta_n(\Gamma)$ is geometrically finite
and every parabolic element of
$\beta_n(\Gamma)$ lies in a rank-two parabolic subgroup,
\item $\{ \beta_n\}$ converges to the identity representation of
$\Gamma$, and
\item if $i_n: M\to M(p_n,q_n)$
denotes the inclusion map, then for each $n$, there exists an
orientation-preserving homeomorphism $\psi_n: {\rm int}
(M(p_n,q_n)) \to \H^3/\beta_n(\Gamma)$ such that
$\beta_n\circ \psi_*$ is conjugate to $(\psi_n)_*\circ (i_n)_*$.
\end{enumerate}
}

\medskip

In order to apply the Hyperbolic Dehn Filling Theorem in our setting we need to
make a couple more topological observations. We choose a meridian-longitude system
for the toroidal boundary component $T$ of $\hat M^\tau_4$ so that the meridian bounds a disk
in $M^\tau_4$ and the longitude bounds an essential annulus $A$ in $\hat M^\tau_4$.
First, one checks that $\hat M^\tau_4(1,n)$ is
homeomorphic to $M^\tau_4$ for all $n$.
(The easiest way to see this is to note that Dehn
twisting $n$ times about $A$ takes the $(1,0)$  curve
on $T$ to the $(1,n)$ curve and observe
that $\hat M^\tau_4(1,0)=M^\tau_4$.) Then one similarly notices that
$(\hat M^\tau_4(1,n), i_n\circ f_\tau)$ is equivalent to $(M^\tau_4,h_\tau)$ in
$\mathcal{A}(M_4)$ for all $n$, where $i_n:\hat M^\tau_4\to M^\tau_4(1,n)$ is the inclusion map
and $h_\tau$ is the homotopy equivalence
defined in example 2 in section \ref{interior}.

Now let $N=\H^3/\Gamma$ be a geometrically finite hyperbolic 3-manifold
admitting an orientation-preserving homeomorphism
$\psi:{\rm int}(\hat M^\tau_4)\to N$. (We further require that all parabolic elements of $\Gamma$
lie in rank two parabolic subgroups.) In the Hyperbolic Dehn Filling Theorem we choose
$(p_n,q_n)=(1,n)$ for all $n$ and let $\beta_n:\Gamma\to {\rm PSL}_ 2({\bf C})$
be the resulting sequence of representations. 

For each $n$, define $\rho_n=\beta_n\circ\psi_*\circ (f_\tau)_*$.  The sequence $\{\rho_n\}$
converges to $\rho=\psi\circ(f_\tau)_*$. Since
$\beta_n\circ \psi_*$ is conjugate to $(\psi_n)_*\circ (i_n)_*$, we
see that $\rho_n$ is conjugate to $(\psi_n)_*\circ (i_n)_*\circ (f_\tau)_*$. 
Since $i_n\circ f_\tau$ is a homotopy equivalence and $\psi_n$
is a homeomorphism, $\rho_n$ is a discrete faithful representation
with image $\pi_1(N_n)=\beta_n(\hat\Gamma)$.
Since $(\hat M^\tau_4(1,n), i_n\circ f_\tau)$ is equivalent to
$(M^\tau_4,h_\tau)$, it follows that $(\hat N_{\rho_n},h_{\rho_n})$
is equivalent to $(M^\tau_4,h_\tau)$ for all $n$. In particular,
$\rho_n$ lies in the component
of ${\rm  int}(AH(M_4))$ associated to $(M^\tau_4,h_\tau)$ for all $n$.
On the other hand, since $f_\tau$ lifts to an embedding, Bonahon'sTameness
Theorem \cite{bonahon} and results of McCullough-Miller-Swarup \cite{MMS} imply
that $N_\rho$ is homeomorphic to ${\rm int}(M_4)$. 
So, $\{\rho_n\}$ is an example of a sequence in $AH(M_4)$ where the
homeomorphism type changes in the limit.

It remains to show that $\rho$ lies in the closure of the
component of ${\rm int}(AH(M_4))$
associated to $(M_4,id)$.  We may accomplish this by modifying the above construction.
We first construct a geometrically finite
hyperbolic 3-manifold $N'$
which admits an orientation-preserving homeomorphism
$\psi':{\rm int}(\hat M_4)\to N'$ (where $\hat M_4$ is
the manifold obtained from $M_4$ by removing an open neighborhood of the core curve
of $V$).
We further require that there exists an embedding $f:M_4\to\hat M_4$ such that
$i_n'\circ f:M_4\to {\rm int}(\hat M_4(1,n))$ is a homotopy equivalence which is
homotopic to an orientation-preserving homeomorphism for all $n$ (where
$i_n':\hat M_4 \to \hat M_4(1,n)$ is the inclusion map). Finally, we require
that $\rho$ is conjugate to $\psi'_*\circ f_*$.
(To construct $N'$  from $M_\rho$, we may normalize so that there is a maximal
abelian subgroup of $\rho(\pi_1(M_4))$ generated by the M\"obius transformation
$z\rightarrow z+1$. If one considers the group $\Gamma'$ generated by
$\rho(\pi_1(M_4))$ and $z\rightarrow z+ri$ for a large enough real value
of $r$, then one may take $N'=\H^3/\Gamma'$.)
One then
applies the Hyperbolic Dehn Filling Theorem just as before to obtain a sequence $\{\rho_n'\}$
of representations lying in the component of ${\rm int}(AH(M_4))$ associated to $(M_4,id)$
and converging to $\rho$. Alternatively, one may simply apply a theorem of Ohshika
\cite{ohshika-boundaries}.

It should be fairly clear that the above construction works for all $n$ and
all $\tau\in\Sigma_n$. In fact, it works for all primitive shuffle
equivalences, see \cite{ACM}.

This construction is also responsible for the phenomenon of self-bumping.
Recall that $M_2=S_3\times I$ where $S_3$ is a closed surface of genus 3.
We illustrate the modifications of the ``wrapping construction''
necessary to produce a self-bumping point in the boundary of $QF(S_3)$.
Let $\hat M_2$ be obtained from $M_2$ by removing an open neighborhood
of the core curve of $V$ and let $(\hat M_2)'$ be the infinite cyclic
cover obtained by gluing alternating copies of $J_1$ and $J_2$ to an
infinite thickened annulus. We construct an orientation-preserving 
embedding  $\hat f:M_2\to (\hat M_2)$
which maps $J_1$ to the ``first copy'' of $J_1$ and maps $J_2$ to the
``second'' copy of $J_2$ (so that there are copies of $J_2$ and $J_1$ which
lie ``between'' the image of $J_1$ and the image of $J_2$).
We then let $f=\pi\circ \hat f$, where $\pi:(\hat M_2)' \to {\hat M}_2$ is
the covering map, and 
apply Thurston's Hyperbolic Dehn Filling Theorem as before, to  produce
a sequence $\{\rho_n\}$ in $QF(S_3)$ 
which converges to a self-bumping point $\rho\in \partial AH(M_2)$.
(The assertion that $\rho$ is a self-bumping point is not obvious and
requires a clever proof, see McMullen \cite{mcmullen} or Bromberg-Holt
\cite{bromberg-holt}).

\section{Relative deformation spaces and the failure of local connectivity}
\label{relative}

There are various naturally defined subsets of $AH(M)$ which have played a prominent
role in the theory of Kleinian groups. Most simply, one may require that certain
elements of $\pi_1(M)$ be mapped to parabolic elements.  In this setting, it is natural
to introduce the language of pared manifolds. We refer the reader to Morgan \cite{Morgan}
and Canary-McCullough \cite{canary-mccullough} for a more extensive discussion of
pared manifolds.

Let $M$ be a compact, orientable, irreducible $3$-manifold with
nonempty boundary which is not a $3$-ball, and let $P\subseteq\partial M$. We say
that $(M,P)$ is a {\em pared} $3$-manifold if 

\begin{enumerate}
\item
Every component of $P$ is an incompressible torus or annulus,
\item
every noncyclic abelian subgroup of $\pi_1(M)$ is conjugate into the
fundamental group of a component of $P$, and
\item
every map $\phi\colon(S^1\times I,S^1\times
\partial I) \to (M,P)$  which induces an injection on fundamental groups
is homotopic, as a map of pairs, to a map $\psi$
such that $\psi(S^1\times I) \subset P$.
\end{enumerate}

We can then define the relative deformation space $AH(M,P)$ to be the
set of (conjugacy classes) of representations $\rho\in AH(M)$ such that if $Q$ is any  component of
$P$ then
$\rho(\pi_1(Q))$ consists of parabolic elements. $AH(M,P)$ lies in a relative
character variety $X(M,P)$. We may define the set $\mathcal{A}(M,P)$ to be
the set of marked, oriented pared manifolds homotopy equivalent to $(M,P)$
up to orientation preserving pared homeomorphisms. Again the interior of
$AH(M,P)$ (within $X(M,P)$) can be identified with
$$\bigcup_{(M',P')\in \mathcal{A}(M,P)} \mathcal{T}(\partial M'-P')/Mod_0(M',P')$$
where $Mod_0(M',P')$ is the group of isotopy classes of  pared self-homeomorphism
of $(M',P')$ which are pared homotopic to the identity
(see Canary-McCullough \cite{canary-mccullough} for more details.)

If $S$ is a compact surface with boundary,
then the interior $QF(S)$ of $AH(S\times I,\partial S\times I)$
is again the space of quasifuchsian groups which are quasiconformally conjugate to
a cofinite area Fuchsian group uniformizing the interior $S^0$ of $S$ and is naturally 
identified with $\mathcal{T}(S^0)\times \mathcal{T}(S^0)$. In the case that $S^0$ is
a once-punctured torus, $AH(S\times I,\partial S\times I)$ is known as the
space of punctured torus groups. Bromberg \cite{bromberg-PT} proved that the space
of punctured torus groups is not locally connected. 

\begin{theorem}{\rm (Bromberg \cite{bromberg-PT})}
\label{bromPT}
The space of punctured torus groups is not locally connected.
\end{theorem}

 Bromberg's proof suggests that many deformation spaces and relative deformation
 spaces are not locally connected. Magid \cite{magid} recently observed that you can
 use Bromberg's result to prove that an infinite family of relative deformation spaces
 is not locally connected. The simplest case of his result
 is the following.

 \begin{theorem}{\rm (Magid \cite{magid})} 
 If $S$ is a closed surface, $M=S\times I$, and $P$ is
 a collection of non-parallel incompressible annuli in $S\times \{ 1\}$
 such that at least one
 component of $S\times \{ 1\}-P$ is homeomorphic to a once-punctured torus,
 then $AH(M,P)$ is not locally connected.
 \end{theorem}

 More generally, Magid's result shows that $AH(M,P)$ fails to be locally
 connected  whenever there is a separating essential
 annulus  $A$ in $M$ with one boundary component lying in a component
 $P_0$ of $P$ such that  the closure $X$  of one  component of $M-A$ 
 is homeomorphic to $S\times I$  where
 $S^0$ is a once-punctured torus 
 (by a homeomorphism taking $\partial S\times I$ to $A$)
 and that $X$ intersects $P$ only in $P_0$.

\medskip

Results of Holt-Souto \cite{holt-souto} and Evans-Holt \cite{evans-holt}
combine to show that the set of self-bumping points, and hence the set of points
where the deformation space is not locally connected, is not dense in the boundary
of the space of punctured torus groups. In particular, they show that there exists some constant
$\epsilon_0>0$ such that if $\rho$ is a punctured torus group and every closed geodesic in
$\partial_cN_\rho$ which is homotopic into a cusp in $\widehat N_\rho$ has length
at least $\epsilon_0$, then the space of punctured torus groups does not self-bump at $\rho$.

Ito \cite{ito-PT} has quite recently given a complete description of the self-bumping points
in the boundary of the space of punctured torus groups. In particular, he shows that all
self-bumping points arise from the construction described in section \ref{key}. Ohshika
\cite{ohshika-exotic} has generalized many of Ito's result to the setting of all quasifuchsian
spaces. Moreover, he is able to show that there is no self-bumping at many points
in the boundary of quasifuchsian space.

\medskip

Another natural and well-studied class of deformation spaces are the Bers slices which
sit inside the quasifuchsian spaces. We recall that
if $S$ is a compact surface, then its associated
space $QF(S)$ of quasifuchsian groups is naturally identified with $\mathcal{T}(S^0)\times\mathcal{T}(S^0)$.
If $\sigma\in\mathcal{T}(S^0)$, then the associated {\em Bers slice} $B_\sigma$ is the set
$\mathcal{T}(S^0)\times \{ \sigma\}$ of quasifuchsian groups, whose bottom conformal
boundary component has conformal structure $\sigma$. Bers \cite{bers-slice} proved
that the closure of any Bers slice is compact, so it is natural to study the topology of the closure
of such a component. Minsky \cite{minsky-PT} proved that if $S^0$ is a once-punctured torus,
then the closure of any Bers slice is homeomorphic to a disk.

\begin{theorem}{\rm (Minsky \cite{minsky-PT})}
If $S^0$ is a once-punctured torus, then the closure in the space of punctured
torus groups of any Bers slice in $QF(S)$ is homeomorphic
to a closed disk.
\end{theorem}

\noindent
{\bf Remark:} Minsky \cite{minsky-PT} also showed that the closure of the Maskit slice is homeomorphic
to a closed disk. Komori \cite{komori} has established the same result for the Earle slice
and Komori and Parkkonen \cite{komori-parkkonen} show that Bers-Maskit slices have
the disk as closure.

\medskip

If the complex dimension of the Bers slice is greater than 1, then very little is known about
the topology of the closure of a Bers slice.  However, Kerckhoff and Thurston \cite{kerckhoff-thurston} showed that, in this
setting, there exist Bers slices $B_\sigma$ and $B_{\sigma'}$ such that the natural homeomorphism
between $B_\sigma$ and $B_{\sigma'}$ does not extend to a homeomorphism of their closures.
They use this to show that there exist Bers slices such that the action of the mapping class
group does not extend continuously to the closure. Kerckhoff and Thurston's results  were the first
indication that the topology of closures of Bers slices must be ``complicated.''

\section{Untouchable points}

In this section, we will describe joint work with Jeff Brock, Ken Bromberg
and Yair Minsky \cite{BBCM} which shows that the topology of $AH(M)$ is well-behaved at ``most''
points in the boundary of $AH(M)$ in many cases.

A point $\rho\in \partial AH(M)$ is said to be {\em untouchable} if there is no bumping or
self-bumping at $\rho$. Notice that $AH(M)$ is locally connected at all untouchable
points.

\begin{theorem}{\rm (\cite{BBCM})}
Suppose that $M$ has incompressible boundary.
If $\rho\in \partial AH(M)$ and $\rho(\pi_1(M))$ contains no parabolic
elements, then $\rho$ is untouchable.
\end{theorem}

If $\partial M$ contains no tori, then such points are generic in $\partial AH(M)$. 

The proof of the non-bumping portion of  this result is a straightforward application of earlier work and we will provide a brief outline of the argument. One first notes that  results of
Thurston \cite{thurston-notes}, Canary \cite{canary-covering} and Anderson-Canary \cite{AC-cores} imply that 
if $\{\rho_i\}\subset \iAH$ converges to $\rho$, then $\{ N_{\rho_i}\}$ converges
geometrically to $N_\rho$ (i.e. larger and larger portions of $N_\rho$ look increasingly
like portions of $N_{\rho_i}$.)
Results of Thurston \cite{thurston-notes}, Canary-Minsky
\cite{canary-minsky}  and Ohshika \cite{ohshika-limits}, then imply that $N_\rho$ is
homeomorphic to $N_{\rho_i}$  (by a homeomorphism in the homotopy class determined
by $\rho\circ\rho_i^{-1}$) for all large $i$. Hence the  (marked) homeomorphism type is
locally constant at $\rho$, so there is no bumping at $\rho$.

The proof that there is no self-bumping at $\rho$ is somewhat more involved. One
begins by showing that if $\{\rho_i\}$ converges to $\rho$, then the 
end invariants of $\{ N_{\rho_i}\}$ converge to those of $N_\rho$. (This is not always the
case, so we are also using our restrictions on  $\rho$ here.) We then consider
any two sequences $\{\rho_i\}$ and $\{\rho_i'\}$ in $\iAH$ converging to $\rho$
and construct paths $\alpha_i:I\to\iAH$ joining $\rho_i$ to $\rho_i'$,
so that the sequence $\{ \alpha_i\} $ converges to the constant path with
image $\rho$.  To establish the convergence of these paths, we
use the same convergence results and arguments as used in the proof
of the Density Theorem as well as the resolution of the Ending Lamination Conjecture.

\medskip

If we allow $N_\rho$ to have cusps, then $\{ N_{\rho_i}\}$ need not converge geometrically
to $N_\rho$ in the above argument, so to rule out bumping we 
must place additional restrictions on $\rho$.
We say that $\rho$ is {\em quasiconformally rigid}, if every component of $\partial_cN_\rho$ is
a thrice-punctured sphere. Notice that this includes the case that $\partial_cN_\rho$ is empty.
(We call such representations quasiconformally rigid, since Sullivan's rigidity theorem
\cite{sullivan-rigid} guarantees that any representation quasiconformally conjugate to
$\rho$ is in fact conformally conjugate.)

\begin{theorem}{\rm (\cite{BBCM})}
If $\rho\in \partial AH(M)$ is quasiconformally rigid, then there is no
bumping at $\rho$.
\end{theorem}

{\bf Remark:} There is a related result of Anderson, Canary and McCullough (Corollary 8.2 in \cite{ACM}) in the
setting where $M$ has incompressible boundary which is much stronger.

\medskip 

We  provide a brief outline of the argument, which is again largely an exercise in applying
previously developed technology.
Suppose that $\{\rho_i\}\subset \iAH$ converges to $\rho$ and
$\{N_{\rho_i}\}$ converges geometrically to $N_\infty$. Then there
exists a covering map $\pi:N_\rho\to N_\infty$.
Results of Anderson-Canary \cite{AC-cores}  and Canary \cite{canary-covering}
imply that there exists a compact core $C$
for $N_\rho$ which embeds in $N_\infty$ (via $\pi$).  For large $i$, one
can pull-back $\pi(C)$ to obtain a compact core $C_i$ for $N_{\rho_i}$
such that $C_i$ is homeomorphic to $C$ (by a homeomorphism in the homotopy class
determined by $\rho_i\circ \rho^{-1}$).  The main result of McCullough-Miller-Swarup \cite{MMS} imply
that $C_i$ is homeomorphic to $N_{\rho_i}$ (by a homeomorphism homotopic to inclusion).
Therefore, $N_{\rho_i}$ is homeomorphic
to $N_\rho$ for all large enough $i$ (again by a homeomorphism in the correct homotopy class).
So, there is no bumping at $\rho$.

\medskip

To rule out self-bumping, we need to further restrict our setting. We recall that a compact,
irreducible manifold is said to be {\em acylindrical} if it does not contain any essential
annuli. (We note that there exist manifolds which are not acylindrical, but do not
contain any {\em primitive} essential annuli.)

\begin{theorem}{\rm (\cite{BBCM})}
If  $M$ is acylindrical  or homeomorphic to $S\times I$ (where $S$ is a closed surface) and 
$\rho\in \partial AH(M)$ is quasiconformally rigid, then there is no self-bumping
at $\rho$. So, $\rho$ is untouchable.
\end{theorem}

The proof of this theorem is much more involved and makes use in a key way of the
techniques of proof of the
Ending Lamination Conjecture, in particular Minsky's a priori bounds \cite{ELC1}, as well as
Thurston's bounded image theorem (see Kent \cite{kent}).

\medskip\noindent
{\bf Remark:}
In the case that $M\cong S\times I$, our results overlap substantially with results
of Ohshika \cite{ohshika-exotic}.  We refer the reader to Ohshika's paper for
the detailed definitions.

\begin{theorem}{\rm (Ohshika \cite{ohshika-exotic})} Suppose that
$S$ is a compact hyperbolic surface and $\rho\in \partial AH(S\times I,\partial S\times I)$.
\begin{enumerate}
\item
If $\partial_cN_\rho$ has one component homeomorphic to $S^0$ and all other
components are thrice-punctured spheres and every rank one cusp of $N_\rho$
(which is not homotopic to a component of $\partial S$)
abuts  a geometrically infinite end, then  there is no self-bumping at $\rho$.
\item
If $\rho$ is quasiconformally rigid and every rank one cusp of $N_\rho$
(which is not homotopic to a component of $\partial S$)
abuts  a geometrically infinite end, then  there is no self-bumping at $\rho$.
\end{enumerate}
\end{theorem}

His techniques are quite different and make deep
use of Soma's work \cite{soma} on geometric limits of quasifuchsian hyperbolic 3-manifolds.
Ohshika and Soma \cite{ohshika-soma} have recently extended Soma's work on geometric limits.

\section{Questions and conjectures}

In this final section, we will collect conjectures and questions about the topology of
deformation spaces of hyperbolic 3-manifolds. As we still understand very little, much remains
to study and we limit ourselves to a few of the more obvious questions.

\subsection{Local connectivity}

Bromberg has conjectured that the failure of local connectivity is a fairly widespread
phenomenon. As a first step, one might hope to show that it holds for all spaces of surface groups.

\begin{conjecture}{\rm (Bromberg \cite{bromberg-PT})}
{\em If $S$ is any compact surface, then $AH(S\times I,\partial S\times I)$ is
not locally connected.}\footnote{Conjecture 10.1 has recently been established by Aaron Magid when $S$ is a
closed orientable surface of genus at least two.}
\end{conjecture}

Bromberg's proof of Theorem \ref{bromPT} makes essential use of the wrapping construction 
described in
section \ref{key}. However, he expects that local connectivity should fail even in settings
where one cannot perform this construction, e.g. for Bers slices. 

\begin{conjecture}{\rm (Bromberg \cite{bromberg-PT})}
{\em If $S$ is any compact hyperbolic surface whose interior $S^0$ is not a once-punctured torus,
thrice-punctured sphere
or 4-times punctured sphere, then the closure of any Bers slice in $QF(S)$ is not locally connected.} 
\end{conjecture}

As Bromberg points out, these two specific conjectures strongly suggest the following
more dramatic conjecture

\begin{conjecture}
{\em If $M$ is any compact hyperbolizable 3-manifold with a non-toroidal boundary component, then
$AH(M)$ is not locally connected.}
\end{conjecture}

It should be pointed out that if $M$ does not contain a primitive essential annulus, then
we don't even know that $M$ is not a manifold, so one might begin with the following
easier conjecture.

\begin{conjecture}
{\em If $M$ is any compact hyperbolizable 3-manifold with a non-toroidal boundary component, then
$AH(M)$ is not a manifold.}
\end{conjecture}

\subsection{The compressible case}

In the case where $M$ has incompressible boundary, Corollary \ref{enum} gives
a complete enumeration of the components of $AH(M)$, but the
situation where $M$ is allowed to have compressible boundary
is still quite mysterious. It is easy to use the constuction in
section \ref{key} to show that various components of $\iAH$ bump.
Bromberg and Holt's result, Theorem 6.2, already applies to
show that if $M$ has compressible boundary, and is not
obtained from one or two thickened tori by adding a single
1-handle, then every component of ${\rm int}(AH(M))$ self-bumps,
since $M$ will contain a primitive essential annulus in
these cases.

One would like to give a complete enumeration of the components of
$AH(M)$ whenever $M$  has compressible boundary. In the incompressible
boundary situation, this is accomplished by showing that the construction
in section \ref{key} is responsible for all bumping phenomena.

\begin{problem}
{\em Give a complete enumeration of the components of
$AH(M)$ in terms of topological data.}
\end{problem}

One might hope to begin by
establishing the following simpler conjecture, which follows from
Theorem \ref{PSE} in the case where $M$ has incompressible boundary.

\begin{conjecture}{\em  $AH(M)$ has finitely many components
if and only if its interior has finitely many components.}
\end{conjecture}

A first step in the proof of this conjecture might be to show that it is not
possible for infinitely many components of $\iAH$ to accumulate at
a single point.

\subsection{The fractal nature of deformation spaces}

Closures of Bers slices
in punctured torus spaces admit natural embeddings in ${\bf C}$. The beautiful pictures
of these embeddings, see Komori-Sugawa-Wada-Yamashita \cite{KSWY}, Yamashita \cite{yamashita},
and Mumford-Series-Wright \cite{indra}, suggest that their boundaries are quite ``fractal''
in nature. One might guess that the boundaries of these
slices have Hausdorff dimension strictly between 1 and 2.

\begin{question} {\em What can one say about the Hausdorff dimension of the boundary
of a Bers slice of a punctured torus?}
\end{question}

Miyachi \cite{miyachi,miyachi2}
has shown that there is a countable dense set of ``cusps'' in the boundary of a Bers
slice of a punctured torus.
Parkkonen \cite{parkkonen} has further analyzed the shapes of these cusps in
the related setting of a Maskit slice. Parkkonen \cite{parkkonen2} has also studied
the shape of Schottky space near certain boundary points.
(Schottky space is ${\rm int}(AH(H_g))$ where $H_g$ is the handlebody of genus $g$.)

Of course, one would like
to study more general classes of deformation spaces in the future.

\subsection{Components of $AH(M)$ and components of the character variety}

It is known that $AH(M)$ can be disconnected, but it is not known
whether or not the different components of $AH(M)$ can lie in different
components of the character variety.

\begin{question} {\em Is it possible that different components of $AH(M)$ lie in different
components of the character variety?}
\end{question}

\end{document}